\documentclass[12pt]{article}
\usepackage{amstext,amsfonts,graphicx,latexsym}
\newtheorem{thm}{Theorem} \newtheorem{lemma}[thm]{Lemma}
\newtheorem{prop}[thm]{Proposition} \newtheorem{cor}[thm]{Corollary}
\newenvironment{pf} {\noindent{\sc Proof. }}{{\hfill
$\Box$}\par\vskip2\parsep} \newenvironment{pfof}[1]
{\par\vskip2\parsep\noindent{\sc Proof of\ #1. }}{{\hfill $\Box$}
\par\vskip2\parsep}

\newcommand{\rd}{{\mathbb R}^d}
\newcommand{\zd}{{\mathbb Z}^d}

\newcommand{\dof}{\bf\boldmath}

\newcommand{\prob}{\mbox{\bf P}}
\newcommand{\expe}{\mbox{\bf E}}
\newcommand{\leb}{{\mathcal L}}

\newcommand{\parend}{{\hfill $\Diamond$}\par\vskip2\parsep}

\newcommand{\esssup}{\mbox{\rm ess sup}}
\newcommand{\as}{a.s{.}\ }

\newcommand{\Z}{\mathbb Z}
\newcommand{\R}{\mathbb R}

\newcommand{\eps}{\epsilon}
\newcommand{\be}{\begin{equation}}
\newcommand{\ee}{\end{equation}}
\newcommand{\C}{\xi} \newcommand{\x}{x} \newcommand{\Y}{y}
\newcommand{\z}{z} \newcommand{\D}{|} \newcommand{\p}{\psi}

\newcommand{\essinf}{\mbox{\rm ess inf}}
\newcommand{\q}{q}
\newcommand{\I}{I--}

\def\qed{\relax\ifmmode\hskip2em \Box\else\unskip\nobreak\hfill
$\Box$\fi}

\newcounter{mycount}

\newenvironment{mylist}{\begin{list}{(\roman{mycount})}%
{\usecounter{mycount}\itemsep 0pt}}{\end{list}}

\title{Tail Bounds for the Stable Marriage of Poisson and Lebesgue}
\author{
Christopher Hoffman\thanks{Funded in part by NSF (USA) grant DMS-0100445
and by MSRI},
Alexander E.
Holroyd\thanks{Funded in part by an NSERC (Canada)
research grant and by MSRI and PIMS}
and Yuval Peres\thanks{Funded in part by NSF (USA) grants DMS-0104073
and DMS-0244479 and by MSRI and CPAM}}

\date{July 15, 2005}

\begin{document}
\maketitle
\renewcommand{\thefootnote}{}
\renewcommand{\thefootnote}{\arabic{footnote}}

\begin{abstract}
Let $\Xi$ be a discrete set in $\rd$. Call the elements of $\Xi$
{\em centers}. The well-known Voronoi tessellation partitions
$\rd$ into polyhedral regions (of varying volumes) by allocating
each site of $\rd$ to the closest center.  Here we study
allocations of $\rd$ to $\Xi$ in which each center attempts to
claim a region of equal volume $\alpha$.

We focus on the case where $\Xi$ arises from a Poisson process of
unit intensity.  It was proved in \cite{hhp} that there is a
unique allocation which is {\em stable} in the sense of the
Gale-Shapley marriage problem.  We study the distance $X$ from a
typical site to its allocated center in the stable allocation.

The model exhibits a phase transition in the appetite $\alpha$. In
the critical case $\alpha=1$ we prove a power law upper bound on
$X$ in dimension $d=1$.  It is an open problem to prove any upper
bound in $d\geq 2$. (Power law lower bounds were proved in \cite{hhp}
for all $d$). In the non-critical cases $\alpha<1$ and $\alpha>1$
we prove exponential upper bounds on $X$.
\end{abstract}

\section{Introduction}

The following model was studied in \cite{hhp}.  Let $d\geq
1$. We call the elements of $\rd$ {\dof sites}. We write $|\cdot
|$ for the Euclidean norm and $\leb$ for Lebesgue measure or {\dof
volume} on $\rd$.  Let $\Xi\subset \rd$ be a
 discrete set. We call the elements of $\Xi$ {\dof centers}. Let
$\alpha\in[0,\infty]$ be a parameter, called the {\dof appetite}. An
{\dof allocation} (of $\rd$ to $\Xi$ with appetite $\alpha$) is a
measurable function
$$\psi:\rd \to \Xi\cup\{\infty,\Delta\}$$
such that $\leb\psi^{-1}(\Delta)=0$, and
$\leb\psi^{-1}(\xi)\leq\alpha$ for all $\xi\in\Xi$.  We call
$\psi^{-1}(\xi)$ the {\dof territory} of the center $\xi$.  We say
that $\xi$ is {\dof sated} if $\leb\psi^{-1}(\xi)=\alpha$, and
{\dof unsated} otherwise.  We say that a site $x$ is {\dof
claimed} if $\psi(x)\in\Xi$, and {\dof unclaimed} if
$\psi(x)=\infty$.

The following definition is an adaptation of that introduced by
Gale and Shapley \cite{gale-shapley}.

\paragraph{Definition of stability.}
Let $\xi$ be a center and let $x$ be a site with
$\psi(x)\notin\{\xi,\Delta\}$.  We say that $x$ {\dof desires}
$\xi$ if
$$|x-\xi| < |x-\psi(x)| \text{ or $x$ is unclaimed.}$$
We say that $\xi$ {\dof covets} $x$ if
$$|x-\xi| < |x'-\xi| \text{ for some $x'\in\psi^{-1}(\xi)$, or $\xi$
is unsated.}$$
We say that a site-center pair $(x,\xi)$ is {\dof
unstable} for the allocation $\psi$ if $x$ desires $\xi$ and $\xi$
covets $x$.  An allocation is {\dof stable} if there are no unstable
pairs.  Note that no stable allocation may have both unclaimed sites
and unsated centers. \parend

Now let $\Pi$ be a translation-invariant, ergodic, simple point
process on $\rd$, with intensity $\lambda\in(0,\infty)$ and law
$\prob$. Our main focus will be on the case when $\Pi$ is a Poisson
process of intensity $\lambda=1$.  The {\dof support} of $\Pi$ is the
random set $[\Pi]=\{z\in\rd:\Pi(\{z\})=1\}$.  We consider stable
allocations of the random set of centers $\Xi=[\Pi]$.

\begin{figure}
\centering

\resizebox{1.6in}{!}{\includegraphics{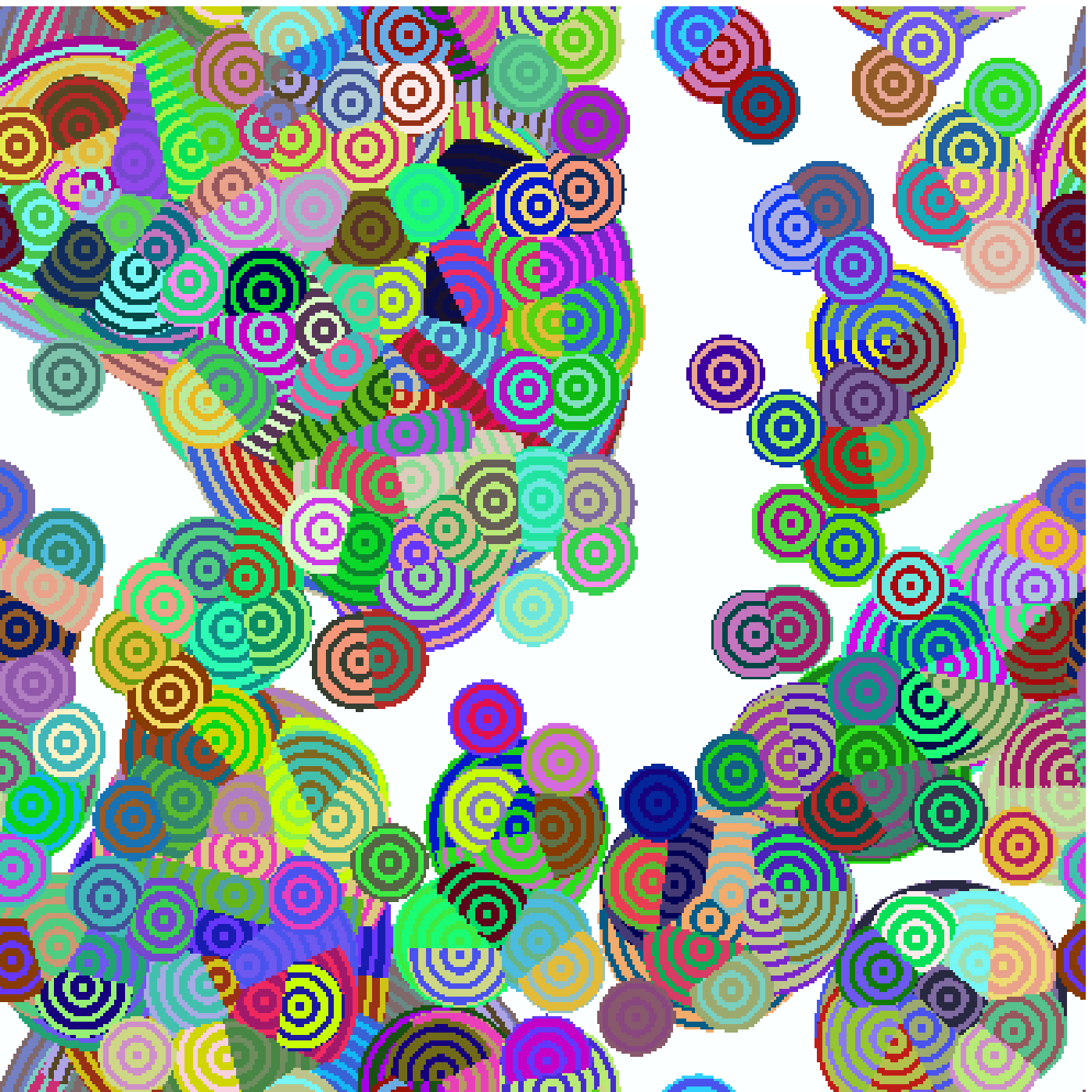}}
\resizebox{1.6in}{!}{\includegraphics{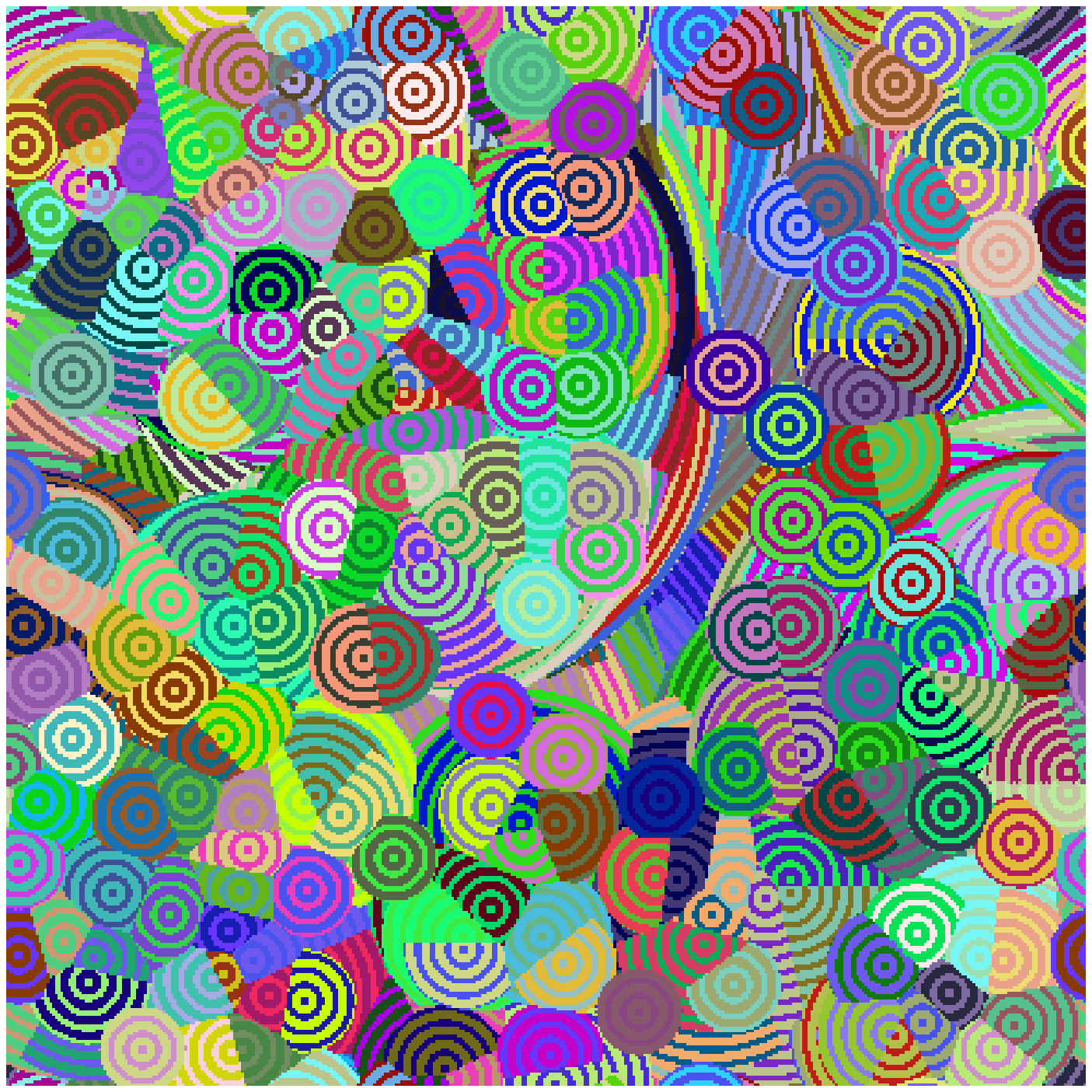}}
\resizebox{1.6in}{!}{\includegraphics{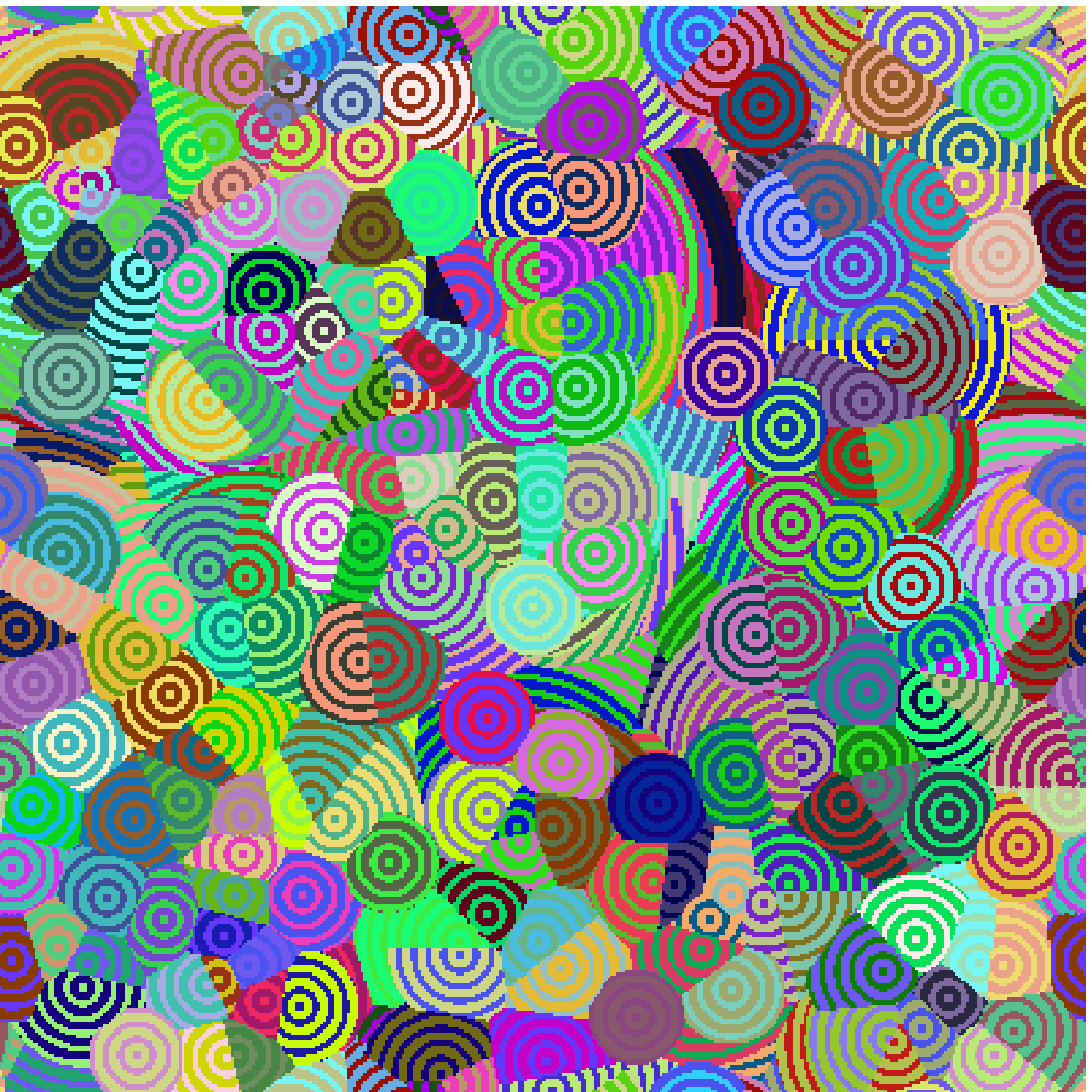}}

\caption{Stable allocations with appetites $\alpha=0.8,\; 1,\;
  1.2$.  The centers are chosen uniformly at random in
  a 2-torus, with one center per unit area.  Each territory is
  represented by concentric annuli in two colors.}
\label{little-sm}
\end{figure}

In \cite{hhp} it was proved that for any ergodic point process
$\Pi$ with intensity $\lambda\in(0,\infty)$ and any appetite
$\alpha\in(0,\infty)$ there is a $\leb$-a.e.\ unique allocation
$\Psi=\Psi_\Pi$ from $\rd$ to $[\Pi]$.  Furthermore we have the
following phase transition phenomenon.
\begin{mylist}
\item If $\lambda\alpha<1$ {\dof (subcritical)} then \as all
centers are sated but there is an infinite volume of unclaimed
sites. \item If $\lambda\alpha=1$ {\dof (critical)} then \as all
centers are sated and $\leb$-a.a. sites are claimed. \item If
$\lambda\alpha>1$ {\dof (supercritical)} then \as not all centers
are sated but $\leb$-a.a. sites are claimed.
\end{mylist}
See Figure \ref{little-sm} for an illustration. For further
information and more pictures see \cite{hhp}.  The critical model
was applied in \cite{h-p-extra} to the construction of certain
shift-couplings.

While the results in (i)--(iii) above suggest the
subcritical / critical / supercritical terminology, the
typical signature of a critical phenomenon in statistical
physics is exponential decay (of correlations, cluster sizes, or large
deviation probabilities) in the subcritical and supercritical regimes,
and sub-exponential decay (usually given by a power law)
at criticality.  We will establish
such a phenomenon for the stable allocation model when
the centers are distributed as a Poisson process.

One natural
quantity to consider is the distance from the origin to its center:
$$X=|\Psi(0)|,$$
where we take $X=\infty$ if $0$ is unclaimed. Another natural
quantity is the {\dof radius} of the territory $\Psi^{-1}(\xi)$:
$$R(\xi)=R_\Psi(\xi)=\esssup_{x\in\Psi^{-1}(\xi)} |\xi-x|.$$ Suppose
$\Pi$ is a Poisson process.
We introduce the point process $\Pi^*$ with law $\prob^*$
obtained from $\Pi$ by adding an extra center at the origin:
\begin{equation}
\label{fake-palm}
[\Pi^*]=[\Pi]\cup\{0\}.
\end{equation}
Define the radius for a typical
center thus:
$$R^*=R_{\Psi_{\Pi^*}}(0).$$
In the subcritical and critical phases, the
conditional law of $X$ given that it is finite is dominated by the law
of $R^*$; see Lemma \ref{xvr} in the remarks below.

\begin{thm}[critical upper bound]
\label{crit-upper}
Let $\Pi$ be a Poisson process with intensity $\lambda=1$.
For $d=1$ and $\alpha=1$ we have $\expe^* (R^*)^{1/18}<\infty.$
\end{thm}

\begin{thm}[non-critical upper bounds] \label{exp-upper}
Let $\Pi$ be a Poisson process \\ with intensity $\lambda=1$.
\begin{mylist}
\item
For all $d$ and $\alpha>1$ we have $\expe e^{c X^d}<\infty$;
\item
For all $d$ and $\alpha<1$ we have $\expe^* e^{c (R^*)^d}<\infty$,
\end{mylist}
for some $c=c(d,\alpha)>0$.
\end{thm}

We shall also prove the following, which answers a question posed
by Lincoln Chayes (personal communication).
\begin{thm}[supercritical rigidity]
\label{rigidity} Let $\Pi$ be a Poisson process with intensity
$\lambda=1$, and consider the stable allocation to the process
$\Pi^*$.  As $\alpha\searrow 1$ we have
$$\prob^*\bigg(0 \text{\rm\ is an unsated center}\bigg)\to 0.$$
\end{thm}

\paragraph{Remarks.}
In the case when $\Pi$ is a Poisson process, the process $\Pi^*$
defined by (\ref{fake-palm}) is the Palm process associated with
$\Pi$; it may be thought of as $\Pi$ conditioned to have a center at
$0$.  The center at $0$ may be thought of as playing the role of a
``typical'' center in the original process $\Pi$.
(The Palm process $\Pi^*$
may also be defined for general point processes, but
(\ref{fake-palm}) is no longer a correct description; see
\cite{hhp},\cite{kallenberg} for more information).

The following simple result relates the random variables $X$ and $R^*$.
\begin{lemma}[site to center comparison]
\label{xvr} Let $\Pi$ be a Poisson process of intensity $\lambda$
and suppose $\lambda\alpha\leq 1$. Then for all $r\in[0,\infty)$
we have
$$\prob(X>r \mid X<\infty) \leq \prob^* (R^*> r).$$
\end{lemma}
Thus, in the subcritical and critical phases, upper bounds for $R^*$
yield corresponding upper bounds for $X$.

In particular, applying Lemma \ref{xvr} to Theorem
\ref{crit-upper} we obtain the following.  Let $\Pi$ be a Poisson
process with intensity $\lambda=1$.
\begin{trivlist}
\item
For $d=1$ and $\alpha=1$ we have $\expe X^{1/18}<\infty.$
\end{trivlist}
It is immediate from Theorem 5(i) of \cite{hhp} that $R^*<\infty$
a.s.\ in all dimensions, but we have been unable to prove any quantitative
upper bound on $R^*$ or $X$ in the critical case in dimensions $d\geq 2$.

The following lower bounds for the critical phase were proved in
\cite{hhp}.  By Lemma \ref{xvr} they imply the analogous lower
bounds for $R^*$.  Let $\Pi$ be a Poisson process with intensity
$\lambda=1$.
\begin{mylist}
\item
For $d=1,2$ and $\alpha=1$ we have $\expe X^{d/2}=\infty$.
\item
For $d\geq 3$ and $\alpha=1$ we have $\expe X^d=\infty$.
\end{mylist}

Applying Lemma \ref{xvr} to Theorem \ref{exp-upper}(ii) we obtain the
following.  Let $\Pi$ be a Poisson process with intensity $\lambda=1$.
\begin{trivlist}
\item
For all $d$ and $\alpha<1$ we have $\expe (e^{c X^d}; X<\infty)<\infty.$
\end{trivlist}
We conjecture that $(R^*)^d$ has a finite exponential moment in
the supercritical case $\alpha>1$ also. It is straightforward to
check that the exponential bounds obtained are tight up to the
value of $c$.  Indeed, denoting the ball
$$B(x,r)=\{y\in\rd: |y-x|<r\},$$
consider the event that $B(0,r)$ contains centers lying
approximately on a densely-packed lattice, while $B(0,2r)\setminus
B(0,r)$ contains no centers. Such an event has probability
decaying at most exponentially in $r^d$ (for any $\alpha$), and it
guarantees that $X>r$ and $R^*>r$.

Our proof of Theorem \ref{exp-upper} does not in general
yield any explicit bound
on the exponential decay constant $c(d,\alpha)$.
However, such a bound is available in each of the following cases:
\begin{mylist}
\item
$d\geq 1$ and $\alpha>2^d$;

\item
$d\geq 1$ and $\alpha<2^{-d}$;

\item
$d=1$ and $\alpha\neq 1$.
\end{mylist}
For the precise statements see Propositions \ref{extreme-alpha} and
\ref{1d-exp}.  The proofs of these results are considerably simpler
than that of Theorem \ref{exp-upper}, and are based on standard
large deviation bounds for the Poisson process.

To what extent are stable allocations robust to changes in the
parameters?  There are several natural ways to formulate such a
question precisely.  We shall prove one such formulation, Theorem
\ref{converge} below, which roughly speaking states that if we change
the set of centers $\Xi$ far away from the origin, then near the
origin the stable allocation $\psi$ changes only on a small
volume.  This result will be a key ingredient in the proofs of
Theorems \ref{exp-upper} and \ref{rigidity}.

In order to state Theorem \ref{converge} precisely,
we need following conventions
(to be used only in Sections \ref{sec:continuity} and \ref{sec:contproof}).
We will work with various
sets of centers, and we want to ensure that they have various
almost sure properties enjoyed by point processes.
We call an allocation $\psi$ to a set of centers $\Xi$ {\dof canonical}
if, for any  $z\in\rd$ and $\zeta\in\Xi\cup\{\infty\}$,
whenever $\leb[B(z,r)\setminus\psi^{-1}(\zeta)]=0$ for some $r>0$ then
$\psi(z)=\zeta$.
We call a set of centers $\Xi$ {\dof benign} if it satisfies
 \begin{mylist}
\item
$\Xi$ has a $\leb$-a.e.\ unique stable allocation, and
\item
$\Xi$ has a unique canonical allocation, which we denote $\psi_\Xi$.
\end{mylist}
By Theorems 1, 3 and 24 of \cite{hhp}, for any ergodic point process $\Pi$
we know that $[\Pi]$ is almost surely a benign set.  (But it
appears hard to describe simple properties of $[\Pi]$ which ensure
that it is benign).  If $\Xi$ is benign then $\psi_\Xi$ has all territories open
and the unclaimed set open.  Furthermore it is the unique minimizer of
the set $\psi^{-1}(\Delta)$ in the class of stable allocations $\psi$
of $\Xi$ with those properties.

For sets of centers $\Xi_1,\Xi_2,\ldots$ and $\Xi$ we write
$\Xi_n\Rightarrow \Xi$ if for any compact $K\subseteq \rd$ there
exists $N$ such that for $n>N$ we have $\Xi_n\cap K=\Xi\cap K$.
For allocations $\psi_1,\psi_2,\ldots$ and $\psi$ we write
$\psi_n\to\psi$ a.e.\ if for $\leb$-a.e.\ $x\in\rd$ we have
$\psi_n(x)\rightarrow\psi(x)$ in the one-point compactification
$\rd\cup\{\infty\}$.
\newpage
\begin{thm}[continuity]
\label{converge} Fix $\alpha$.  Let $\Xi_1,\Xi_2,\ldots$ and $\Xi$
be benign sets of centers, and write $\psi_n=\psi_{\Xi_n}$ and
$\psi=\psi_\Xi$ for their canonical allocations.  If
$$\Xi_n\Rightarrow \Xi$$
then
$$\psi_n \rightarrow \psi \qquad
\text{a.e.}$$
\end{thm}

We shall refer extensively to results from the companion article
\cite{hhp}.  We adopt the convention that ``Theorem I--x'' refers to
Theorem x of \cite{hhp}.

\section{Site to Center Comparison}

\begin{pfof}{Lemma \ref{xvr}}
Our proof applies in the more general context when $\Pi$ is an
ergodic point process of intensity $\lambda\in(0,\infty)$, and
$\Pi^*$ is the Palm process (see \cite{hhp},\cite{kallenberg}
for more details).
First note that by Theorem \I 4, $\lambda\alpha\leq 1$ implies
that all centers are sated a.s.

We claim that for any $r\in(0,\infty)$,
\begin{equation}
\label{equal-in-d}
\prob\bigg(R(\Psi_\Pi(0))>r \;\bigg|\; X<\infty\bigg)=\prob^*(R^*>r),
\end{equation}
so $R(\Psi_\Pi(0))$ conditioned on $0$ being claimed is equal in
distribution to $R^*$.  Once this is proved, the result follows
because clearly
$$X\leq R(\Psi_\Pi(0)), \qquad \prob\text{-a.s.}$$

We shall use the mass-transport principle (Lemma \I 17).  For
$z\in\zd$ let $Q_z=[0,1)^d+z\subset\rd$, and define
$$m(u,v)=\expe\leb\bigg\{x\in Q_u:\Psi_\Pi(x)\in Q_v,
 \;R(\Psi_\Pi(x))>r\bigg\}.$$
Using Fubini's Theorem and translation invariance we have
\begin{eqnarray*}
\sum_{v\in\zd} m(0,v)
&=&\expe\leb\bigg\{x\in Q_0: \Psi_\Pi(x)\neq \infty
\text{ and }R(\Psi_\Pi(x))>r\bigg\} \\
&=& \prob\bigg(X<\infty\text{ and }R(\Psi_\Pi(0))>r\bigg).
\end{eqnarray*}
On the other hand, since all centers are sated, and by a standard
property of the Palm process,
$$
\sum_{u\in\zd} m(u,0)
= \alpha\expe\#(\xi\in[\Pi]\cap Q_0: R(\xi)>r)=\alpha\lambda\prob^*(R^*>r).
$$
Lemma \I 17 states that $\sum_{v\in\zd} m(0,v)=\sum_{u\in\zd}
m(u,0)$, and by Proposition \I 20 we have
$\prob(X<\infty)=\alpha\lambda$, so (\ref{equal-in-d}) follows.
\end{pfof}

\section{One Dimensional Critical Bound}
\label{sec:1d}

In this section we deduce Theorem \ref{crit-upper} from a more
general result.  Let $\Pi$ be a stationary renewal process, and
let $\Pi^*$ be its Palm version. Write the support
$[\Pi^*]=\{\xi_j: j \in \Z \}$, where $(\xi_j)$ is an increasing
sequence. Thus $(\xi_j)$ is a two-sided random walk, with
$\xi_0=0$. We assume that the i.i.d.\ increments $\xi_j-\xi_{j-1}$
have mean 1 and finite variance $\sigma^2$. In the (critical)
stable allocation with $\alpha=1$, our goal is to prove a power
law tail bound for $R^*=R_{\Psi_{\Pi^*}}(0)$.

\begin{thm} \label{thm:renew}
With the assumptions above, there exists a constant $C<\infty$
that depends on the law of $\xi_j-\xi_{j-1}$, such that for all
$r>1$,
$$
\prob^*(R^* >r) \le C r^{-1/17.6} \,.
$$
\end{thm}

\begin{pfof}{Theorem \ref{crit-upper}}
This is immediate from Theorem \ref{thm:renew}.
\end{pfof}

\begin{figure}
\centering
\includegraphics{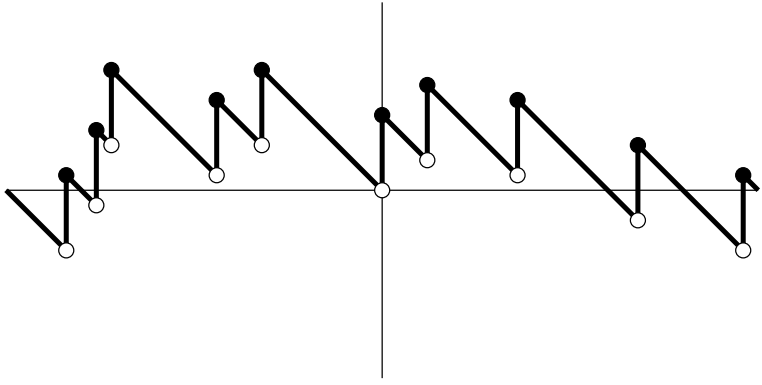}
\caption{The function $F$ jumps up by $1$ at each center, and has
slope $-1$ elsewhere.} \label{function}
\end{figure}

\medskip
Given $\Xi \subseteq \R$, write $N(s,t]=\#(\Xi\cap(s,t])$.
 We introduce the function $F:\R\to\R$ defined by
\begin{eqnarray*}
F(0)&=&1; \\
F(t)-F(s)&=&N(s,t]-(t-s) \qquad\text{for }s<t.
\end{eqnarray*}
See Figure \ref{function} for an illustration.
 Note that for all $\xi \in \Xi$, we have
$F(\xi+)=F(\xi)=F(\xi-)+1$.

We will prove that if $F$ has certain properties then $R^*$
cannot be too large.  On the other hand, we can analyze the
behavior of $F$ using the technology of random walks.

\begin{prop}[measure-preserving map] \label{cor.pres}
Let $\Xi\subseteq\R$ be a discrete set of centers and let $\psi$
be a stable allocation to $\Xi$ with appetite $\alpha=1$ in which
all centers are sated.  For each center $\xi$, the restriction of
$F$ to $\psi^{-1}(\xi)$ is a measure-preserving map into
$[F(\xi-),F(\xi)]$ (where on both sides, the measure is $\leb$).
\end{prop}
We will prove the above proposition from the following.
\begin{lemma} \label{l.measure}
Let $\Xi\subseteq\R$ be a discrete set of centers and let $\psi$
be a stable allocation to $\Xi$ with appetite $\alpha=1$.
\begin{mylist}

\item Suppose that $t>\xi$ satisfies $\psi(t)=\xi$ and
$\leb\Bigl(\psi^{-1}(\xi) \cap [\xi,t] \Bigr)=\beta$. Then
$F(t)=F(\xi)-\beta$.

\item Suppose that $s<\xi$ satisfies $\psi(s)=\xi$ and
$\leb\Bigl(\psi^{-1}(\xi) \cap [s,\xi] \Bigr)=\gamma$. Then
$F(s)=F(\xi-)+\gamma$.
\end{mylist}
\end{lemma}
\begin{pf}
By symmetry, it suffices to prove (i). \newline \noindent{\bf Case
A:} Suppose $F(t)>F(\xi)-\beta$. Then $N(\xi,t]>t-\xi-\beta$, so
there is a center $\eta \in (\xi,t]$ and a set of positive length
of sites $s \notin (\xi,t]$ with $\psi(s)=\eta$. If $s>t$ or
$\eta-s>t-\eta$ then $(t,\eta)$ is an unstable pair, so we must
have $s \in [2\eta-t, \xi)$ for a.e.\ such $s$. But then $(s,\xi)$
is an unstable pair.

\noindent{\bf Case B:} Suppose that $F(t)<F(\xi)-\beta$. Then
$N(\xi,t]<t-\xi-\beta$, so there exists a center
$\eta\notin[\xi,t]$ and a set of positive length of sites $s \in
(\xi,t)$ with $\psi(s)=\eta$. If $\eta<\xi$ or $\eta-s>s-\xi$,
then $s,\xi$ is an unstable pair, so we must have $\eta \in [t,
2s-\xi)$ for a.e.\ such $s$. But then $t,\eta$ is an unstable pair.
\end{pf}

\begin{pfof}{Proposition \ref{cor.pres}}
Denote  $\Gamma_+=\psi^{-1}(\xi)\cap(\xi,\infty)$
and $\gamma_+=\leb(\Gamma_+)$.
Part (i) of Lemma \ref{l.measure} implies that
 the restriction of $F$ to $\Gamma_+$ is a monotone decreasing,
measure-preserving map into $[F(\xi)-\gamma_+,F(\xi)]$.
Part (ii) of Lemma \ref{l.measure} implies that
 the restriction of $F$ to
$\Gamma_-=\psi^{-1}(\xi)\cap(-\infty,\xi)$ is a monotone
decreasing, measure-preserving map into
$[F(\xi-),F(\xi-)+\gamma_-]$, where $\gamma_-=\leb(\Gamma_-)$.
Since $\gamma_-+\gamma_+=\alpha=1$, this completes the proof.
\end{pfof}

\begin{lemma} \label{l.match}
Under the assumptions of Proposition \ref{cor.pres}, suppose that
$0\in\Xi$ and that $x>0$ is such that $F(t)<0$ for $t \in (x,2x)$.
Then $R(0)  \le x$.
\end{lemma}

\begin{pf}
Let $r \le x$ be maximal such that $F(r-)=0$. Denote $D=\{y \in
(0,r) \, : \, 0<F(y)<1 \}$. Consider two cases.

\noindent{\bf Case I:} There exists a site $t \in D$ such that
$\psi(t)<0$ or $\psi(t)>2x$. In this case, since there is a center
at 0, stability of the pair $(0,t)$ implies that $0$ must be sated
by distance $r$.

\noindent{\bf Case II:} $\psi(D) \subseteq [0,2x]$. Then Lemma
\ref{l.measure} implies that $\psi(D) \subseteq [0,r]$.  Since
$F(r)=F(0-)=0$, the upcrossings and downcrossings of $[0,1]$ by
$F$ from $0$ to $r$ must equalize, i.e.
\begin{equation} \label{cross}
\sum_{\xi \in \Xi\cap [0,r]} \Bigl[(F(\xi)\wedge 1)-(F(\xi-)\vee
0)\Bigr] =\leb(D) \,.
\end{equation}
By Proposition \ref{cor.pres}, for every center $\xi \in [0,r]$ we
have $$
 (F(\xi)\wedge 1)-(F(\xi-)\vee 0) \ge
\leb(D\cap\psi^{-1}\xi) \, ,
$$
and the identity in (\ref{cross}) implies that for each $\xi$ the
above inequality must be an equality. In particular, for $\xi=0$
this shows that $R(0) \le r \le x$.
\end{pf}

The following random walk lemma will provide the tail estimate needed
to prove Theorem \ref{thm:renew}.
\begin{lemma}[random walk estimate] \label{l.walk}
Let $\{X_j\}_{j \ge 1}$ be i.i.d.\ random variables, with mean
zero and variance $\sigma^2<\infty$. Suppose $X_j \ge -1$ a.s.
Write $S_k=\sum_{j=1}^k X_j$ and denote
\begin{equation} \label{defAm}
A_m:=\Bigl\{S_j >1 \; \mbox{ \rm for all }  j \in [2^{3m-1},
2^{3m}) \Bigr\}  \,.
\end{equation}
Then $$\prob\bigg(\bigcap_{k=1}^m A_k^c\bigg) \le C_1\theta^m$$
 for some $\theta<8^{-1/35.2}$ and
$C_1<\infty$.
\end{lemma}

\begin{pf}
Write $D_m=\bigcap_{k=1}^m A_k^c$.
 It suffices to show that for all $m$ sufficiently large,
\begin{equation} \label{cond}
\prob\Bigl(A_{m+1}^c \, \Big| D_m \Bigr) \le \theta\,.
\end{equation}
Fix $m>1$, and denote $M:=2^{3m}$. Consider the event
\begin{equation} \label{defGm}
G_m:=\Bigl\{\exists j\in [M,2M] : S_j  \ge -1 \Bigr\} \,.
\end{equation}

We will first show that
\begin{equation} \label{cond2}
\prob\Bigl(G_m\, \Big| D_m \Bigr) \ge 1-2^{-1/2}+o(1) \;
  \mbox{ \rm as } m \to \infty\,.
\end{equation}
It clearly suffices to show this when $D_m$ on the right-hand side
is replaced by $D_m\cap\{S_M<0\}$. To do so, let $\tau$ be the
largest integer $j \le M$ such that $S_j \ge 0$. Denote by
$\tau^*$ the index of the last maximum for the walk in $(\tau,M]$,
so that $S_i \le S_{\tau^*}$ for $i \in (\tau,\tau^*]$, and
$S_{\tau^*} > S_j$ for $j \in (\tau^*,M]$. Note that
$S_{\tau^*}\geq -1$, since all $X_j\geq-1$.  We will derive
(\ref{cond2}) from the uniform estimate
\begin{equation} \label{cond3}
\prob\Bigl(G_m\, \Big| D_m, \tau^*, S_{\tau^*} \Bigr) \ge
1-2^{-1/2}+o(1) \,.
\end{equation}
Observe that conditional on $\tau^*=\ell$, the sequence
$\{S_{\ell+i}-S_\ell\}_{i=0}^{\infty}$ has the same law as  the
sequence $\{S_i \}_{i=0}^{\infty}$ conditioned to stay negative
for the interval $i\in[1,M-\ell]$, and this also applies when we
condition further on $D_m$ and on the value of $S_\ell$. By
\cite{feller-ii} Chapter\ XII formula (8.8),
$$
\prob\Bigl(S_i<0 \mbox{ \rm for }  i \in [1, k) \Bigr)=
(c+o(1))k^{-{1/2}}
$$
as $k \to \infty$, and furthermore the probability is non-zero for all
$k\geq 1$.  Therefore,
\begin{eqnarray*} \label{grow} &
\prob\Bigl(S_{\ell+i} < S_\ell \;\;\forall i \in (M-\ell,2M-\ell]
\;\Big|\; S_{\ell+i} < S_\ell \;\;\forall i \in [1,M-\ell] \Bigr) & \\
\ &\displaystyle\leq
\frac{(c+o(1))(2M-\ell)^{-1/2}}{(c+o(1))(M-\ell-1)^{-1/2}} & \\ &\leq
2^{-1/2}+o(1) & \end{eqnarray*} uniformly in $\ell$ as $M=2^{3m} \to
\infty$. This proves (\ref{cond3}) and hence (\ref{cond2}).

Next, we show that
\begin{equation} \label{grow2}
\prob\Bigl(A_{m+1} \, \Big|\, G_m\cap D_m \Bigr) \ge
(1/\pi)\arcsin 3^{-1/2}-o(1) \,.
\end{equation}
Indeed by the strong Markov property, it suffices to show that
\begin{equation} \label{grow3}
\prob\Bigl(A_{m+1} \, \Big| S_j \ge -1 \Bigr) \ge (1/\pi)\arcsin
3^{-1/2}-o(1) \,.
\end{equation}
holds uniformly in $j \in [M,2M]$. By Donsker's Theorem
(\cite{kallenberg} Theorem 14.9) this, in turn, is a consequence
of the following inequality for standard Brownian motion
$B(\cdot)$. For all $u \in [1,2]$ and $b \ge 0$,
\begin{eqnarray*}
\prob\Bigl(\min_{4 \le t \le 8} B(t) \ge 0 \,\Big|\, B(u)=b  \Bigr)
&\geq&
\prob\Bigl(\min_{4 \le t \le 8} B(t) \ge 0 \,\Big|\, B(2)=0 \Bigr) \\
& = & \prob\Bigl(\min_{1/3 \le t \le 1} B(t) \geq 0\Bigr) \\
&=& (1/\pi)\arcsin 3^{-1/2}-o(1).
\end{eqnarray*}
The latter follows from the arcsine law for the last zero of
Brownian motion on an interval (\cite{kallenberg} Theorem 13.16).

In conclusion, we obtain  (\ref{cond}), with any $\theta$ such
that $1-\theta< (1-2^{-1/2}) (1/\pi)$ $\arcsin 3^{-1/2}$. This is
valid if $\theta>8^{-1/35.2}$, proving the lemma.
\end{pf}

\begin{pfof}{Theorem \ref{thm:renew}}
Let $X_j=\xi_j-\xi_{j-1}-1$ so that $S_j=\xi_j-j$. On the event $A_m$
defined in Lemma \ref{l.walk}, we have $\xi_j >j+1$ for all $j \in
[M/2,M)$, where $M=2^{3m}$.  Therefore, on $A_m$ we have $F(t)<0$
for all $t \in [M/2,M)$, whence $R(0)<M/2$ by Lemma \ref{l.match}.
Therefore
$$
\{R(0) \ge M/2\} \subseteq \bigcap_{k=1}^m A_k^c \, .
$$
  So far, we
have only considered the centers on the positive axis, and our
estimates hold uniformly over the positions of centers on the
negative axis. By considering the symmetrical events on the
negative axis, we obtain
$$
\prob^*[R(0) \ge M/2] \le \prob^*(D_m)^2 \le C_0 8^{-m/17.6} \le C_0
M^{-1/17.6} \, .
$$
Given any $r>1$, we can choose $m$ maximal so that $M/2=2^{3m-1}
\le r$. Since $C_0 M^{-1/17.6} \le C r^{-1/17.6}$ for a suitable
$C$, the theorem follows.
\end{pfof}

\section{Explicit Exponential Bounds}

In this section we prove exponential upper bounds involving explicit
constants in several cases.  Denote $\q(x)=[x-1-\log(x)]/x$ so that
$\q(x)>0$ for all positive $x \ne 1$.  Write $\omega_d$ for the volume
of the unit ball in $\R^d$.

\begin{prop} [explicit bounds for extreme $\alpha$]
\label{extreme-alpha} Let $\Pi$ be a Poisson \\ process on $\R^d$
with intensity $1$.
\begin{mylist}
\item
For any $\alpha>2^d$ we have
$\expe e^{c X^d}<\infty$  provided $c<\omega_d \q(\alpha/2^d)$.
\item
For any $\alpha<2^{-d}$ we have
$\expe^* e^{c (R^*)^d}<\infty$  provided
$c<2^d\omega_d \q(\alpha 2^d)$.
\end{mylist}
\end{prop}

\begin{prop} [explicit bounds for $d=1$]
\label{1d-exp} Let $\Pi$ be a Poisson process on $\R$ with
intensity $1$. For any $\alpha \ne 1$, we have
$$\text{
$\expe( e^{c X} ; X<\infty)<\infty$  and $\expe^* e^{c R^*}<\infty$
provided $c<\q(\alpha)$.
}$$
\end{prop}

(In fact, the proofs of Propositions \ref{extreme-alpha} and
\ref{1d-exp} give explicit upper bounds on the tail probabilities
$\prob(X>r)$ and $\prob^*(R^*>r)$.)

\begin{pfof}{Proposition \ref{extreme-alpha}}
We first note a standard large deviation estimate.
If $Z$ is a Poisson
 random variable with mean $\gamma$ we have:
\begin{eqnarray}
\label{potop} \prob(Z \ge b) \le e^{-\gamma \q(\gamma/b)} \;
 &\quad\text{ \rm for}& b>\gamma ; \\
\label{pobot} \prob(Z\le a) \le e^{-\gamma \q(\gamma/a)} \;
 &\quad\text{ \rm for}& 0<a<\gamma \, .
\end{eqnarray}
Indeed, (\ref{potop}) follows from setting $s=b/\gamma$ in
$ s^b \prob(Z \ge b) \le \expe(s^Z)=e^{\gamma(s-1)} \,,$
and (\ref{pobot}) follows similarly from
$(a/\gamma)^a \prob(Z \le a) \le \expe((a/\gamma)^Z)$.
See e.g. \cite{kallenberg} Chapter 27.

For (i), fix $\alpha>2^d$ and let $Z$ be the number of centers
in $[\Pi] \cap B(0,r)$. Then $Z$ is Poisson with mean $\omega_d r^d$.
 On the event that $Z>\omega_d r^d 2^d/\alpha$,
there must be at least one some center $\xi$ in $[\Pi] \cap B(0,r)$
which is not sated within $B(0,2r)$. Stability of the pair $(0,\xi)$
then implies that $0$ must be allocated to some center
no farther than $\xi$, whence $X \le |\xi| <r$.
Thus $\prob(X>r) \le \prob(Z \le \omega_d r^d 2^d/\alpha)$;
an application of (\ref{pobot}) completes the proof.

For (ii), fix $\alpha<2^{-d}$ and let $Z'$ be the number of centers
in $[\Pi^*] \cap B(0,2r)$. Then $Z'-1$ is a Poisson with mean
$\omega_d 2^d r^d$. On the event that $Z'<\omega_d r^d /\alpha$,
there must be (a positive volume of) sites $x$ in $B(0,r)$
that are not allocated to any center in $[\Pi] \cap B(0,2r)$.
Stability of such a site $x$ and the center $0$
implies that $0$ must be sated within the closed ball
$\overline{B(0,|x|)}$, whence $R^* \le |x| <r$.
Thus $\prob^*(R^*>r) \le \prob(Z'-1 \ge \omega_d r^d/\alpha-1)$;
an application of (\ref{potop}) completes the proof.
\end{pfof}

In order to prove Proposition \ref{1d-exp}, it will be convenient to
work with $\alpha=1$ and $\lambda\neq 1$ and then rescale.  Recall the
definition of the function $F$ from Section \ref{sec:1d}.  The
following states that sites are allocated to centers on the same level
of $F$.

\begin{lemma}
\label{same-level} Let $\Xi\subseteq\R$ be a discrete set of
centers and let $\psi$ be a stable allocation to $\Xi$ with
appetite $\alpha=1$.  If $\psi(x)=\xi\in\Xi$ then
$F(x)\in[F(\xi-),F(\xi)]$.
\end{lemma}

\begin{pf}
The result is immediate from Lemma \ref{l.measure}, since for any
interval $I$ we have $\leb(\psi^{-1}(\xi)\cap I)\in[0,1]$.
\end{pf}

\begin{pfof}{Proposition \ref{1d-exp}}
We start by noting the following standard large deviation estimates.
If $\Pi$ is a Poisson process with intensity $\lambda$ on $\R$, then
for any $r,a\geq 0$ we have:
\begin{eqnarray}
\label{pp-dev-1}
\prob(\exists t\geq r : \Pi(0,t]\leq t+a)\leq
 \lambda^a e^{-q(\lambda)\lambda r}
&\quad\text{for}& \lambda>1; \\
\label{pp-dev-2}
\prob(\exists t\geq r : \Pi(0,t]\geq t-a)\leq
\lambda^{-a}e^{- q(\lambda)\lambda r}
&\quad\text{for}& \lambda<1.
\end{eqnarray}
To prove the above facts, consider the martingale
$$M(t):=e^{t(\lambda-1)}\lambda^{-\Pi(0,t]}.$$
If $\lambda>1$,
consider the stopping time $\tau=\inf\{t\geq r: \Pi(0,t]\leq t+a\}$.
On the event $\tau=t$, where $t<\infty$, we have $M(\tau)\geq
e^{t(\lambda-1)}\lambda^{-(t+a)}= \lambda^{-a}e^{q(\lambda)\lambda t}$.
Hence applying the optional stopping theorem to $\tau\wedge N$ yields
$1=\expe M(\tau\wedge N)\geq \prob(\tau<N)\lambda^{-a} e^{
q(\lambda) \lambda r}$, and taking $N\to\infty$ yields (\ref{pp-dev-1}).  For
(\ref{pp-dev-2}) we apply similar reasoning to $\tau'=\inf\{t\geq r:
\Pi(0,t]\geq t-a\}$.

Now we prove exponential bounds on $X$ and $R^*$ in the
case when $d=1$, $\alpha=1$ and $\lambda\neq 1$; then we will rescale $\R$.

Firstly, let $\Xi=[\Pi]$.  By Lemma \ref{same-level}, on the event
that $r<X<\infty$ there exists some center $\xi\in [\Pi]\setminus[-r,r]$
with $F(\xi)\in[0,1]$.  Recalling the definition of $F$, taking $t=|\xi|$
and using (\ref{pp-dev-1}),(\ref{pp-dev-2})
we therefore obtain
\begin{eqnarray*}
\lefteqn{\prob(r<X<\infty)} \\
&\leq&
\prob(\exists t>r : \Pi(0,t]\in[t-1,t])
+ \prob(\exists t>r : \Pi(-t,0]\in[t-1,t]) \\
&\leq& 2(1\vee \lambda^{-1}) e^{- q(\lambda)\lambda r}.
\end{eqnarray*}

Secondly, let $\Xi=[\Pi^*]$.  By Lemma \ref{same-level}, on the event
that $R^*>r$ there exists $x\in \R \setminus[-r,r]$
with $F(x)\in[0,1]$,  so
we obtain similarly
\begin{eqnarray*}
\lefteqn{
\prob^*(R^*>r)} \\
& \leq & \prob^*(\exists t>r : \Pi^*(0,t]\in[t-1,t])
+  \prob^*(\exists t>r : \Pi^*(-t,0]\in[t-1,t]) \\
&\leq& 2(1\vee \lambda^{-2}) e^{- q(\lambda)\lambda r}.
\end{eqnarray*}

Finally, rescaling $\R$ by a factor of $\lambda$ changes the intensity
to $1$ and the appetite to $\lambda$, while scaling $X$ and $R^*$ by a
factor of $\lambda$.  Thus we obtain the desired results.
\end{pfof}

\section{Continuity}
\label{sec:continuity}

Recall the continuity result, Theorem \ref{converge}, stated in the
introduction.  In this section we deduce some consequences which will
be used in the proofs of Theorems \ref{exp-upper} and \ref{rigidity}.
The proof of Theorem \ref{converge} is deferred until the end of the
article.

We shall apply Theorem \ref{converge} as follows.  Roughly
speaking, given an almost sure local property of stable
allocations, we may find some large box such that with high
probability the property holds throughout the box, whatever the
configuration of centers outside. More precisely, we apply this to
the notions of replete sets and decisive sets as described below.

In what follows we take $\alpha=1$, and take $\Pi$ to be a Poisson
process of intensity $\lambda$, with associated probability measure
and expectation $\prob_\lambda,\;\expe_\lambda$.  Lemma \ref{triv} and
Corollary \ref{box} below apply to the critical and subcritical
models; that is to $\lambda\leq 1$.  The critical case will be used to
prove Theorem \ref{rigidity} and the subcritical case will be used to
prove Theorem \ref{exp-upper}(ii).

Recall that $\psi_\Xi$ denotes the canonical allocation of the
benign set of centres $\Xi$.
Given a benign set $\Xi\subseteq\rd$ and a measurable $A
\subseteq \R^d$, let $\Xi'_A$ be a random set of centers which is
the union of $\Xi\cap A$ and a Poisson process of intensity
$\lambda$ in $\R^d \setminus A$.  Write $\mu_{\lambda,A}$ for the
law of $\Xi'_A$.  For $\xi\in \Xi$, say that $A$ is $\Xi$-{\dof
replete} for $\xi$ if for every $\lambda\in(0,\infty)$ we have for
$\mu_{\lambda,A}$-a.e.\ $\Xi'_A$ that $\Xi'_A$ is benign and
$\leb\left(\p^{-1}_{\Xi'_A}(\xi)\cap A \right)=1$.
(That is, if $\xi$ is sated
within $A$ whatever happens outside $A$).

Define the box $Q(L)=[-L,L)^d$.
\begin{lemma} \label{triv}
Let $\alpha=1$ and let $\Pi$ be a Poisson process of intensity
$\lambda\leq 1$.  Let $G$ be the event that for every $\xi\in[\Pi]$
there exists $L<\infty$ such that $\xi+Q(L)$ is $[\Pi]$-replete for
$\xi$. Then $\prob_\lambda(G)=1$.
\end{lemma}

\begin{cor}[replete boxes] \label{box}
Let $\alpha=1$ and let $\Pi$ be a Poisson process of intensity
$\lambda\leq 1$.  For any $\epsilon >0$ there exists $M$ such that
$$
\expe_\lambda \#\bigg\{ \xi\in [\Pi]\cap Q(M) : Q(M) \text{\rm\ is not
$[\Pi]$-replete for }\xi\bigg\}<\epsilon (2M)^d.
$$
\end{cor}

Now given benign $\Xi$, we say that a measurable set $A\subseteq\rd$ is
$\Xi$-{\dof decisive} for a site $x$ if for every
$\lambda\in(0,\infty)$ we have $\mu_{\lambda,A}$-a.s.\ that
$\Xi'_A$ is benign and
$\psi_{\Xi'_A}(x)=\psi_\Xi(x)$. (That is, if $\psi(x)$ can be
determined by looking only at $\Xi\cap A$).  Note that if $A$ is
$\Xi$-decisive for $x$ then $\psi_\Xi(x)$ cannot be a center
outside $A$.

The supercritical case below will be used to prove Theorem
\ref{exp-upper}(i).

\begin{lemma}
\label{decisive} Let $\alpha=1$ and let $\Pi$ be a Poisson process
of intensity $\lambda\geq 1$.  Then $\prob_\lambda$-a.s.\ there exists
$L<\infty$ such that $Q(L)$ is $[\Pi]$-decisive for $0$.
\end{lemma}

\begin{cor}[decisive boxes]
\label{decisive-box} Let $\alpha=1$ and let $\Pi$ be a Poisson
process of intensity $\lambda\geq 1$.  For any $\epsilon>0$ there
exists $M<\infty$ such that
$$\expe_\lambda
\leb\bigg[x\in Q(M):Q(M) \text{\rm\ is not $[\Pi]$-decisive for
}x\bigg]<\eps (2M)^d.
$$
\end{cor}

Next we turn to the proofs of the four results above.

\begin{lemma}
\label{extra-bit} Suppose $\Xi_n\Rightarrow \Xi$ and
$\psi_n\to\psi$ a.e.\ are as in Theorem \ref{converge}.  If there is a
set $A$ of positive volume such that every $z\in A$ desires $\xi$
under $\psi$, then for $n$ sufficiently large, $\xi$ is sated in
$\psi_n$, and
$$\limsup_{n\to\infty} R_{\psi_n}(\xi)\leq \essinf_{z \in
A}|z-\xi|\;\;(<\infty).$$
\end{lemma}

\begin{pf}
As the set $A$ has positive volume,  Theorem \ref{converge} implies
that there exists $\z \in A$ such that $\p_n(\z)\to \p(\z)$. Thus
for $n$ sufficiently large, $z$ desires $\xi$ under $\psi_n$.
By stability $\xi$ does not covet $z$, and the result follows.
\end{pf}

\begin{pfof}{Lemma \ref{triv}}
On $G^c$, there exists a center $\xi\in[\Pi]$
  such that for each $L$ there is a benign set of centers $\Xi_L$
  agreeing with $[\Pi]$ on $\xi+Q(L)$ and satisfying
\begin{equation}
\label{funny}
\leb\left(\p^{-1}_{\Xi_L}(\xi)\cap (\xi+Q(L)) \right)<1
\end{equation}
for each $L$.  Suppose that $\prob_1(G^c)>0$ and write
$\psi_L=\psi_{\Xi_L}$. Since $\Xi_L\Rightarrow [\Pi]$, Theorem
\ref{converge} implies that $\psi_L\to\psi_{[\Pi]}$ a.e.  Furthermore,
 Lemma \ref{extra-bit} applies to $\xi$ (by Theorem \I 4(i)
 if $\lambda<1$ or Theorem \I 6(i) if $\lambda=1$), so almost surely
 for $L$ sufficiently
large $\xi$ is sated in each $\psi_L$, and the radii $R_{\psi_L}(\xi)$
are bounded as $L \to \infty$.  This contradicts (\ref{funny}).  We
conclude that $\prob_\lambda(G^c)=0$.
\end{pfof}

\begin{pfof}{Corollary \ref{box}}
For $A \subseteq \R^d$, let $\Pi^L(A)$ denote the number of
$\xi\in [\Pi]\cap A$ such that $ \xi+Q(L)$ is not $[\Pi]$-replete
for $\xi$.
 Lemma \ref{triv} and
the monotone convergence theorem imply that $ \expe_\lambda
\Pi^L(Q(1)) \to 0$ as $L \to \infty$. Thus we can choose an
$L<\infty$ so that the translation-invariant point process
$\Pi^L$ has intensity less than $\eps/2$. Observe that for $M>L$
and $\xi \in [\Pi]\cap Q(M-L)$, if
 $Q(M)$ is not  $[\Pi]$-replete for $\xi$, then $\xi \in [\Pi^L]$.
Therefore
\begin{eqnarray*}
&\expe_{\lambda} \#\{ \xi\in [\Pi]\cap Q(M) :  Q(M) \text{ is not
$[\Pi]$-replete for }\xi\}& \\
 &< (\eps/2) (2M-2L)^d+ (2M)^d-(2M-2L)^d, &
\end{eqnarray*}
which is smaller than $\epsilon (2M)^d$ if $M$ is sufficiently large.
\end{pfof}

In order to prove Lemma \ref{decisive} we need the following
enhancement of Theorem \ref{converge}, in which we (partially) specify
the set on which a.e.\ convergence occurs.  The proof is deferred until
the end of the article.
\begin{prop}
\label{converge-explicit}
 Suppose $\Xi_n\Rightarrow \Xi$ and $\psi_n\to\psi$ a.e.\ are as
in Theorem \ref{converge}. If $\p(z)=\xi\in\Xi$ and $z$ is not
equidistant from any two centers of $\Xi$ then $\p_n(z)\to \xi.$
\end{prop}

\begin{pfof}{Lemma \ref{decisive}}
Since $\lambda\geq 1$, by Theorem \I 4, $0$ is claimed a.s. And
a.s.\ $0$ is not equidistant from any two centers.

Now, on the complement of the event in the lemma, for each $L$
there exists a benign $\Xi_L$ agreeing with $[\Pi]$ on $Q(L)$ such
that $\psi_{\Xi_L}(0)\neq\psi_{[\Pi]}(0)$.

But by Proposition \ref{converge-explicit}, when the all the
events mentioned above occur we have
$\psi_{\Xi_L}(0)\to\psi_{[\Pi]}(0)$ as $L\to\infty$, a
contradiction.
\end{pfof}

\begin{pfof}{Corollary \ref{decisive-box}}
Fix $\epsilon>0$.  Let $U^L$ be the (random) set of sites $x$ for
which $Q(L)+x$ is not $[\Pi]$-decisive.  Then the process $U^L$ is
translation-invariant in law, and by Lemma \ref{decisive}, we may
fix $L$ large enough so that it has intensity less that $\eps/2$.
Now if $M$ is sufficiently large then
$$\expe_\lambda\leb[U^L\cap Q(M)]<(\eps/2)(2M)^d,$$
whence
\begin{eqnarray*}
&  \expe_\lambda \leb[x\in Q(M):Q(M) \text{\rm\ is not
$[\Pi]$-decisive for }x] & \\
&<(\eps/2) (2M)^d + (2M)^d - (2M-2L)^d, &
\end{eqnarray*}
which is less than $\eps (2M)^d$ if $M$ is sufficiently large.
\end{pfof}

\section{Supercritical Rigidity}

In this section we prove Theorem \ref{rigidity}.

\begin{lemma}[coupling]
\label{coupling} For any set $A\subseteq\rd $ of finite volume and
any $\delta
>0$ there exist $\lambda>1$ and a coupling $(\Pi_1,\Pi_\lambda)$
of two Poisson processes of respective
intensities $1,\lambda$, such that
$$\expe\bigg[\Pi_\lambda(A);
\Pi_\lambda\neq\Pi_1 \text{\rm\ on }A\bigg]<\delta.$$
\end{lemma}

\begin{pf}
We take $\Pi_\lambda=\Pi_1+\Pi_\beta$ where $\Pi_1,\Pi_\beta$ are
independent Poisson processes of intensities $1,\beta$ with
$1+\beta=\lambda$.  Then we have
\begin{eqnarray*}
\expe[\Pi_\lambda(A);\Pi_\lambda\neq\Pi_1 \text{ on }A] &=&
\expe[\Pi_1(A)+\Pi_\beta(A);\Pi_\beta(A)>0] \\
&=&
\expe[\Pi_1(A)]\prob[\Pi_\beta(A)>0]+\expe[\Pi_\beta(A)] \\
&=& (\leb A)(1-e^{-\beta\leb A})+\beta\leb A \\
&\to & 0 \qquad\text{ as }\beta\to 0.
\end{eqnarray*}
\end{pf}

\begin{pfof}{Theorem \ref{rigidity}}
By rescaling $\rd$, the required result is equivalent to the same
limiting statement as $\lambda\searrow 1$ with $\alpha=1$.

Given any $\epsilon>0$, choose $M$ by Corollary \ref{box} so that,
writing $\Pi_1=\Pi$,
\begin{equation} \label{newline}
\expe_{1} \#\bigg\{ \xi\in[\Pi_1]\cap Q(M): Q(M) \text{ is not
$[\Pi_1]$-replete for }\xi\bigg\}<\epsilon (2M)^d.
\end{equation}
 Then choose $\lambda$ and a coupling $(\Pi_1,\Pi_\lambda)$
 by Lemma \ref{coupling} so that
\begin{equation} \label{newline2}
\expe\bigg(\#([\Pi_\lambda]\cap Q(M)); \Pi_\lambda\neq\Pi_1 \text{
on }Q(M)\bigg)<\epsilon (2M)^d.
\end{equation}
Note that if $\Xi=\Xi'$ on $Q(M)$ then $Q(M)$ is $\Xi'$-replete
for a center if and only if it is $\Xi$-replete. Therefore, using
(\ref{newline}) and (\ref{newline2}),
\begin{eqnarray*}
\lefteqn{ \expe_{\lambda} \#\bigg\{ \xi\in[\Pi_\lambda]\cap Q(M):
Q(M) \text{ is not $[\Pi_\lambda]$-replete for }\xi\bigg\} }
 \\
 & \leq  &
\expe\bigg( \#\big\{ \xi\in[\Pi_1]\cap Q(M): Q(M)
\text{ is not $[\Pi_1]$-replete for }\xi\big\}; \\
& & \qquad\qquad\qquad\qquad\qquad
  \Pi_1=\Pi_\lambda \text{ on }Q(M)\bigg) \\
& &+ \expe\bigg( \#( [\Pi_\lambda]\cap Q(M) );
  \Pi_1\neq \Pi_\lambda \text{ on }Q(M)\bigg) \\
&\leq & 2 \epsilon (2M)^d.
\end{eqnarray*}
So in particular
$$\expe_\lambda\#(\xi\in [\Pi_\lambda]\cap Q(M): \xi \text { is
unsated})\leq 2\epsilon (2M)^d,$$
and therefore since $\Pi^*$ is
the Palm process,
$$\prob^*_\lambda(0\text{ is unsated})\leq 2\epsilon.$$
\end{pfof}

\section{Supercritical Bound}

In this section we prove Theorem \ref{exp-upper}(i).  Let
$\alpha=1$ and let $\Pi$ be a Poisson process of rate $\lambda$
with law $\prob_\lambda$.

\begin{thm}
\label{exp-decay-super}
Let $\alpha=1$ and let $\Pi$ be a Poisson process of intensity $\lambda$.
For any $\lambda >1$ there exist $C,c\in(0,\infty)$ such
that for all $r>0$,
$$\prob_{\lambda}(X>r)<Ce^{-cr^d}.$$
\end{thm}

\begin{pfof}{Theorem \ref{exp-upper}(i)}
By rescaling $\rd$, the required result is equivalent to the same
statement with $\alpha=1$ and $\lambda>1$, and this is immediate
from Theorem \ref{exp-decay-super}.
\end{pfof}

\begin{pfof}{Theorem \ref{exp-decay-super}}
First observe that if
\begin{equation}
\label{newstar} \text{there exists }\xi\in[\Pi]\cap B(0,r) \text{
with } \leb[\Psi^{-1}(\xi)\cap B(0,2r)]<1
\end{equation}
then
$$X\leq r.$$
This is because $\xi$ must covet some $z\notin B(0,2r)$, so
$|\xi-z|>r$; but $|0-\xi|<r$, so $(0,\xi)$ would be unstable if
$X>r$.

So it is enough to show that the probability that (\ref{newstar})
fails decays exponentially in $r^d$.  Given $\lambda$ let
$$\epsilon=\frac{\lambda-1}{10 \cdot 2^d}\wedge 1,$$
and let $M=M(\lambda,\epsilon)$ be as in Corollary
\ref{decisive-box}.  Note that $\epsilon$ and $M$ do not depend on
$r$.

Now for any $r>0$ we tile the shell $B(0,2r)\setminus B(0,r)$ with
disjoint copies of the box $Q(M)$.  Recall that $Q(M)=[-M,M)^d$.
For $\z\in\zd$ write
$Q_z=Q(M)+2M z$, and define the random variable
$$Y_z=\leb(x\in Q_z: Q_z \text{ is not $[\Pi]$-decisive for $x$}).$$
Let
$$I=I(r)=\{z\in\zd: Q_z\subseteq B(0,2r)\setminus B(0,r)\}$$
be the index set of the boxes lying entirely in the shell, and let
$$S=S(r)=[B(0,2r)\setminus B(0,r)]\setminus \bigcup_{z\in I} Q_z$$
be the remainder of the shell.

Observe that if $r$ is sufficiently large then
\begin{equation}
\label{kind-of-F} \leb S<\epsilon \leb B(0,r).
\end{equation}
Also consider the events
\begin{eqnarray*}
E&=&\bigg\{\Pi( B(0,r))>(\lambda-\epsilon)\leb B(0,r)\bigg\};\\
G&=&\bigg\{\sum_{z\in I}Y_z < 4\epsilon 2^d \leb B(0,r)\bigg\}.
\end{eqnarray*}
We claim that if $E$ and $G$ occur and (\ref{kind-of-F}) holds
then (\ref{newstar}) is satisfied.
To verify this claim, note that given those assumptions,
\begin{eqnarray*}
\leb[x\in B(0,2r):\Psi(x)\in B(0,r)] &\leq &
\sum_{z\in I} Y_z + \leb B(0,r)+\leb S \\
&\leq & (4\epsilon 2^d +1 +\epsilon)\leb B(0,r) \\
 &< & (\lambda-\epsilon)\leb B(0,r) \\
 &< & \Pi(B(0,r))
\end{eqnarray*}
(Here the third inequality holds because by the choice of
$\epsilon$ we have $4\eps 2^d+2\eps < 10\eps 2^d\leq \lambda-1$).
Then recalling that $\alpha=1$ we see that (\ref{newstar}) must
indeed hold.

Finally, we must show that  $\prob(E^C),\prob(G^C)$ each decay at
least exponentially in $r^d$ as $r\to\infty$.  For $E^C$ this is a
standard large deviations bound since $\Pi(B(0,r))$ is Poisson
with mean $\lambda\leb B(0,r)=\Theta(r^d)$ as $r\to\infty$.
Turning to $G^C$, note that the random variables $(Y_z)_{z\in I}$
are i.i.d.\ with mean less than $\epsilon (2M)^d$ by Corollary
\ref{decisive-box}. We have $\# I=\Theta(r^d)$, while
$$(\# I)(2M)^d \leq \leb B(0,2r)\leq 2\cdot 2^d \leb B(0,r),$$
and hence
$$\expe_\lambda\bigg(\sum_{z\in I} Y_z\bigg)\leq (\# I) \epsilon (2M)^d
\leq 2 \epsilon 2^d \leb B(0,r).$$
 Furthermore, each random variable $Y_z$ is bounded by $(2M)^d$.
Therefore by the Chernoff bound (\cite{kallenberg} Corollary
27.4), $\prob(G^C)$ decays exponentially in $r^d$.
\end{pfof}

\section{Subcritical Bound}

In this section we prove Theorem \ref{exp-upper}(ii), via the
following.

\begin{thm} \label{exp-decay-sub}
Let $\alpha=1$ and let $\Pi$ be a Poisson process of intensity $\lambda$.
For any $\lambda <1$ there exist $C,c>0$ such that for all $r>0$,
$$\prob_{\lambda}(\exists \C \in [\Pi]\cap B(0,1) \mbox{ such that
}R(\C)>r)<Ce^{-cr^d}.$$
\end{thm}

\begin{pfof}{Theorem \ref{exp-upper}(ii)}
First note that by rescaling $\rd$, it suffices to prove the same
statement for $\alpha=1$ and $\lambda<1$.  Let $C,c$ be as in
Theorem \ref{exp-decay-sub}. Let $Y$ be the number of centers
$\xi\in [\Pi]\cap B(0,1)$ with $R(\xi)>r$, and note that by a
standard property of the Palm process, $\expe (Y)=\lambda \leb
B(0,1)\prob^*(R^*>r)$, so it is enough to prove that $\expe(Y)$
decays exponentially in $r^d$. Let $u=e^{c r^d/2}$. Then note that
\begin{eqnarray*}
\expe(Y)&=&\expe(Y; 0<Y\leq u)+\expe(Y; Y>u) \\
&\leq & u\prob(Y>0)+\expe[\Pi(B(0,1));\Pi(B(0,1))>u].
\end{eqnarray*}
From Theorem \ref{exp-decay-sub} we have $\prob(Y>0)\leq Ce^{-c
r^d}$, while the second term is bounded above by
$\expe(\Pi(B(0,1))^2)/u$.  Thus both terms decay exponentially in
$r^d$, hence so does \expe(Y).
\end{pfof}

\begin{pfof}{Theorem \ref{exp-decay-sub}}
Fix $\lambda<1$. First observe that if
\begin{equation}
\label{star5}
\text{there exists } y\in B(0,r) \text{ with } \Psi(y)\notin B(0,2r+1)
\end{equation}
then
$$R(\xi)<r+1 \text{ for all } \xi\in [\Pi]\cap B(0,1).$$
This is because otherwise we would have $|y-\xi|<r+1$ and
$|y-\Psi(y)|>r+1$, and so $(y,\xi)$ would be unstable.

So it is enough to show that the probability (\ref{star5}) fails
decays exponentially in $r^d$.  Let
$$\epsilon=\frac{1-\lambda}{10\cdot 2^d},$$
 and let $M=M(\lambda,\epsilon)$ be as in Corollary \ref{box}.
 Note that $\epsilon$ and $M$ do not depend on $r$.

Now for any $r>0$ we tile the shell $B(0,2r+1)\setminus B(0,r)$
with disjoint copies of the box $Q(M)$.  For $\z\in\zd$ write
$Q_z=Q(M)+2M z$, and define the random variable
$$W_z=\#(\xi\in[\Pi]\cap Q_z: Q_z \text{ is not $[\Pi]$-replete
 for $\xi$}).$$
Let
$$I=I(r)=\{z\in\zd: Q_z\subseteq B(0,2r+1)\setminus B(0,r)\}$$
be the index set of the boxes lying entirely in the shell, and let
$$S=S(r)=[B(0,2r+1)\setminus B(0,r)]\setminus \bigcup_{z\in I} Q_z$$
be the remainder of the shell.

Consider the events
\begin{eqnarray*}
E&=&\bigg\{\Pi( B(0,r))<(\lambda+\epsilon)\leb B(0,r)\bigg\};\\
F&=&\bigg\{\Pi(S)<\epsilon \leb B(0,r)\bigg\}; \\
G&=&\bigg\{\sum_{z\in I}W_z < 4\epsilon 2^d \leb B(0,r)\bigg\}.
\end{eqnarray*}
We claim that if $E$, $F$ and $G$ all occur
then (\ref{star5}) is satisfied.
To verify this claim, recall that $\alpha=1$, so that on $E$ we have
\begin{eqnarray*}
\leb\{\Y \in B(0,r) :\p(\Y) \not\in B(0,r)
\}&\geq &
\leb B(0,r)-(\lambda+\epsilon)\leb B(0,r) \\
&\geq & 9\epsilon 2^d \leb B(0,r),
\end{eqnarray*}
(The second inequality holds because $9\epsilon 2^d+\epsilon\leq
10\epsilon 2^d = 1-\lambda$ by the choice of $\epsilon$).  On $F$
we clearly have
$$\leb\{\Y \in B(0,r) :\p(\Y) \in S \}<\epsilon\leb B(0,r),$$
while on $G$, by the definition of replete we have
$$\leb\{\Y \in B(0,r) :\p(\Y) \in \cup_{z\in I}Q_z \}<4\epsilon
2^d\leb B(0,r).$$ Therefore since $B(0,2r+1)=S \cup B(0,r) \cup
\bigcup_{z\in I}Q_z$, on $E\cap F\cap G$ we have
$$\leb\{\Y\in B(0,r) :\p(\Y) \not\in B(0,2r+1) \}
    \geq (9\epsilon 2^d-4\epsilon 2^d-\epsilon) \leb B(0,r)>0,$$
establishing the claim.

Finally, we must show that  $\prob(E^C),\prob(F^C),\prob(G^C)$
each decay at least exponentially in $r^d$ as $r\to\infty$.  For
$E^C$ this is a standard large deviations bound since
$\#([\Pi]\cap B(0,r))$ is Poisson with mean $\lambda\leb
B(0,r)=\Theta(r^d)$. For $F^C$ it also follows from the standard
large deviations bound on noting that $\leb S< (\epsilon/2)\leb
B(0,r)$ for $r$ sufficiently large. Turning to $G^C$, note that
the random variables $(W_z)_{z\in I}$ are i.i.d.\ with mean less
than $\epsilon (2M)^d$ by Corollary \ref{box}.  We have $\#
I=\Theta(r^d)$, while
$$(\# I)(2M)^d \leq \leb B(0,2r+1)\leq 2\cdot 2^d \leb B(0,r),$$
and hence
$$\expe_\lambda\bigg(\sum_{z\in I} W_z\bigg)\leq (\# I) \epsilon
(2M)^d \leq 2 \epsilon 2^d \leb B(0,r).$$ Furthermore, we have
$W_z\leq \Pi(Q_z)$ so each random variable $W_z$ has exponentially
decaying tails.  Therefore by the Chernoff bound
(\cite{kallenberg} Corollary 27.4), $\prob(G^C)$ decays
exponentially in $r^d$.
\end{pfof}

\section{Proofs of Continuity Results}
\label{sec:contproof}

\begin{pfof}{Theorem \ref{converge}}
We can find a countable dense set $X\subseteq\rd$ such
that $\p_{n}(\x)\neq \Delta$ for each $\x \in X$ and for all
$n$. We can choose a subsequence $(n_j)$ such that $\p_{n_j} (\x)$
converges in the compact space $\Xi \cup \infty$
 for all $\x \in X$. We define the map
$\p_\infty$ by
$$\p_{\infty}(\z)=\lim_{j\to\infty} \p_{n_j}(\z)$$
for all $\z$ where the limit exists. Thus $\p_\infty$ exists on
$X$ and perhaps elsewhere.

We define $$\widetilde{R}_{\infty}(\xi)= \sup\bigg\{|\x-\xi
    |:x\in X \text{ and }\p_{\infty}(\x)=\xi \bigg\}.$$
Let
\begin{eqnarray*}
Z&=&\bigcup_{\xi} \bigg\{w\in \rd:|w-\xi|=\widetilde{R}_\infty(\xi)\bigg\} \\
&\cup&\bigcup_{\xi\neq\xi'} \bigg\{w\in\rd: |w-\xi|=|w-\xi'|\bigg\},
\end{eqnarray*}
where the first and second unions are over all centers and all pairs of
 centers in $\Xi$ respectively.  And let
$$D=\bigcup_{n} \psi_n^{-1}(\Delta).$$
The sets $Z$ and $D$ are $\leb$-null a.s.

For $z\in\rd$ let
$$S(\z)=\{\C \in \Xi \cup \infty:\ \exists x_1,x_2,\ldots \in X
\text{ such that  $\x_j\to \z$ and } \p_{\infty} (\x_j) \to \C\}.$$
 By the
compactness of $\Xi \cup \infty$, for any $\z$ the set $S(\z)$ is
not empty.  We claim the following.
\paragraph{Claim.}
\begin{equation}
\label{one-element}
\text{If $z\notin Z\cup D$ then $\psi_\infty(z)$ exists
and $S(z)=\{\psi_\infty(z)\}$.}
\end{equation}
To prove this, we take $z\notin Z\cup D$ and consider two cases.

\paragraph{Case I.}
Suppose that $\C \in S(z) \setminus\{\infty\}$. Since $\C \in S(\z)$
we have that
$$\D\C-z| \leq \widetilde{R}_{\infty}(\xi),$$
and as $\z \notin Z$ we deduce
$$\D\C-z| < \widetilde{R}_{\infty}(\xi).$$
Hence we can pick $\x\in\psi_\infty^{-1}(\xi)\cap X$
 such that $\D\x -\C|>\D\z
-\C|$.   Since $\p_\infty (\x)$ exists there is $N$ such that we
have that $\p_{n_j} (\x)=\C$ for all $n_j>N$, so $\xi$ covets $z$
under $\psi_{n_j}$.  Since $(z,\xi)$ is stable for
$\psi_{n_j}$ we deduce that $\D \p_{n_j}(\z)-z| \leq \D\z -\C|$
for all $n_j>N$.

Label $\Xi\cap \overline{B(z,|z-\xi|)}=\{ \C_1
,\C_2,\dots ,\C_\ell\}$ in such a way that
$$\D\z -\C_1|<\D\z -\C_2|< \cdots < \D\z -\C_\ell|.$$
This is possible as $\z \notin Z$.  (Note that $\xi=\xi_\ell$).
Furthermore since $\z \notin Z$ there exists $r>0$ such that:
\begin{enumerate}
\item $r<\min_i \bigg|\widetilde{R}_{\infty}(\xi_i)-\D\z -\C_i|
\bigg|$, and \item $\D y -\C_1|<\D y -\C_2|< \cdots <\D y -\C_\ell|$
for all $y\in B(z,r)$.
\item $|y-\eta|>|y-\xi|$  for all
$\eta\in\Xi\setminus \{\xi_1,\ldots, \xi_\ell\}$ and all $y\in B(z,r)$.
\end{enumerate}
We will show that for
$\leb$-a.e.\ $y\in B(z,r)$ we have $\psi_{n_j}(y)\to\xi$.

Let $L=\min\{i:\widetilde{R}_\infty(\xi_i)>|z-\xi_i|\}$.  We first show
that $\psi_{n_j}(y)\to\xi_L$ for all $y\in B(z,r)\setminus D$.  By
the definition of $r$ there exists $w \in X$ and $N$ such that
for all $n_j>N$
$$\psi_{n_j}(w)=\xi_L$$ and
\begin{equation}
\label{closer}
|z-\xi_L|+r<|w-\xi_L|.
\end{equation}
For $n_j>N$ and
for every $y\in B(z,r)$ with $\psi_{n_j}(y)\neq\Delta$,
from the stability of $(y,\xi_L)$ under $\psi_{n_j}$, and by (\ref{closer})
we have
$$|y -\psi_{n_j}(y)|<|y-\xi_L|\leq |w-\xi_L|.$$
Therefore for all $y\in B(z,r)\setminus D$ we have
$$\psi_{n_j}(y)=\xi_{i_j}$$
for $n_j>N$, where $i_j=i_j(y) \leq L$. Our next task is to show
that in fact $i_j<L$ is impossible for $j$ sufficiently large.

Suppose on the contrary that there exists $I<L$ and a subsequence
$(n_{j_k})$ and sites $(y_{j_k})$ such that for all $k$
$$\psi_{n_{j_k}}(y_{j_k})=\xi_I$$
and
$$|y_{j_k}-\xi_I|>|z-\xi_I|-r.$$
Then there exists $u \in X \cap B(z,r)$ such that for all $k$
$$|u-\xi_I|<|y_{j_k}-\xi_I|.$$
Since $u \in X$, the sequence $\psi_{n_{j_k}}(u)$
 converges to some $\xi_i$.  By
stability of $(u,\xi_I)$ under $\psi_{n_{j_k}}$ and by the choice
of $r$ we must have $i \leq I<L$. Thus
$$\widetilde{R}_{\infty}(\xi_i) \geq |u-\xi_i|>|z-\xi_i|-r.$$  By
the choice of $r$ the previous line implies
$\widetilde{R}_{\infty}(\xi_i)\geq |z-\xi_i|+r.$ This contradicts
the definition of $L$, so there is no $I<L$ as described.

We have shown that for all $y \in B(\z,r)\setminus D$ the
sequence $\p_{n_j}(y)$ converges to the same center
$\xi_L$, and since $\xi\in S(z)$, this center must be
 $\C$.  Since $z\notin D$ we have that $\psi_\infty(z)=\xi$.
Hence we have proved
claim (\ref{one-element}) in Case I.

\paragraph{Case II.}
Suppose $S(\z)\cap\Xi=\emptyset$; then $S(z)=\{\infty\}$, and we
want to show that $\p_{\infty}(\z)=\infty$.  We work by
contradiction.  Suppose there exists $\C \in \Xi$ and a
subsequence
 $(n_{j_k})$ such that
$\p_{n_{j_k}}(\z)\to \C$.
 Then there exists $r>0$ such that for all
$\x \in X \cap B(\z,r)$
$$|\p_{\infty} (\x)-\x|>|\C -z |.$$  As $\z
\notin \Xi$ we may further choose $\x \in X \cap B(\z,r)$ such that
Then there exists $j_k$ such that $(\x,\C)$ is an unstable pair for
$\p_{n_{j_k}}$.  Hence we have proved claim (\ref{one-element})
in Case II also.
\vspace{2\parsep}

We have proved that $\p_\infty$ is defined almost everywhere.  It
is straightforward to show that if $\Xi_n\Rightarrow \Xi$ and
$\psi_n\to\psi_\infty$ then $\psi_\infty$ is a stable allocation
to $\Xi$ (the main step is to show that $\psi^{-1}(\xi)=\liminf
\psi_n^{-1}(\xi)=\limsup \psi_n^{-1}(\xi)$ a.e.).  Since $\Xi$ is
benign it has an a.e.\ unique stable allocation, so $\psi_\infty$
must agree with $\psi$ a.e. Thus we have $\psi_{n_j}\to\psi$ a.e.

Finally we prove convergence of the entire sequence.  We claim
that for all $z\in\rd\setminus D$ satisfying
$$\psi(z)\neq \Delta$$
and
\begin{equation}
\label{non-equid-z}
|z-\xi|\neq |z-\xi'| \text{ for all } \xi\neq\xi'\in\Xi,
\end{equation}
we have $\psi_n(z)\to\psi(z)$.  Suppose this does not hold for some
$z$ where $\psi(z)=\zeta\in\Xi\cup\{\infty\}$ say.  Then there exists
$(n_j)$ such that
\begin{equation}
\label{wrong-limit}
\psi_{n_j}(z)\neq \zeta \text{ for all } n_j.
\end{equation}
Also since $\psi$ is a canonical allocation we have
\begin{equation}
\label{good-nbrd} \psi(y)=\zeta \text{ for all $y$ in a
neighborhood of $z$}.
\end{equation}
We will deduce a contradiction.

First suppose $\zeta\neq\infty$.  As before, using
(\ref{non-equid-z}) we can label $\Xi\cap\overline{B(z,|z-\xi|)}$
$=\{ \C_1 ,\C_2,\dots ,\C_\ell\}$ with $\xi_\ell=\zeta$ and choose
$r>0$ so that for all $y\in B(z,r)$ we have
\begin{equation}
\label{in-ball} \D y -\C_1|<\D y -\C_2|< \cdots < \D y -\C_\ell|.
\end{equation}

By $(\ref{good-nbrd})$ and the subsequential convergence proved
earlier there exist $x_0,\ldots ,x_\ell\in B(z,r)$ and a subsequence
$(n_{j_k})$ such that
\begin{mylist}
\item
$\psi_{n_{j_k}}(x_i)=\zeta$ for all ${j_k}$ and $i=0,\ldots,\ell$,
\item
$|x_i-\xi_i|<|z-\xi_i|$ for $i=1,\ldots,\ell$, and
\item
$|x_0-\zeta|>|z-\zeta|$.
\end{mylist}
By (i) $(i=0)$, (iii) and stability we have for all $j_k$,
$$|\psi_{n_{j_k}}(z)-z|\leq |\zeta-z|.$$
By (i) $(i>0)$, (ii), (\ref{in-ball}) and
stability we have
$$\psi_{n_{j_k}}(z)\neq\xi_i \text{ for any }i=1,\ldots, \ell-1.$$
Thus for all $j_k$ we have $\psi_{n_{j_k}}(z)=\zeta$ which contradicts
(\ref{wrong-limit}).

Finally suppose $\zeta =\infty$. If $\p_n(z)$ does not converge to
$\infty$ then there exists a subsequence $\p_{n_j}$ and center
$\xi$ such that $\p_{n_j}(z)= \xi$ for all $j$.  By
$(\ref{good-nbrd})$ and the subsequential convergence proved
earlier there exist $x$ and a further subsequence $n_{j_k}$ such
that
$$|x-\xi|<|z-\xi| $$
and
$$\p_{n_{j_k}}(x) \to \infty. $$
Thus for $k$ large enough we have that
$$|\p_{n_{j_k}}(x)-x|>|\xi-x|.$$
 By stability for these $k$ we have that
 $$\p_{n_{j_k}}(z)\neq \xi.$$
 This is a contradiction.
\end{pfof}

\begin{pfof}{Proposition \ref{converge-explicit}}
Assume that $\p(z)=\xi$.
Label $\Xi \cap \overline{B(z,|z-\xi|)}$ as
$\xi_1,\dots,\xi_\ell=\xi$ such that
\be \label{linestar6}
|z-\xi_1|<|z-\xi_2|<\dots<|z-\xi_\ell|.
\ee
As $\p^{-1}(\xi)$ is open there exists $r>0$ such that
\be \label{linestar5}B(z,r) \subset \p^{-1}(\xi).\ee
As $|z-\xi'|\neq |z-\xi''|$ for all $\xi' \neq \xi''$, we can
choose $r$ such that (\ref{linestar5}) is satisfied and for all $y
\in B(z,r)$
\be \label{linestar4}
|y-\xi_1|<|y-\xi_2|<\dots<|y-\xi_\ell|.
\ee

As $\p_n$ converges a.e.\ we can find $y_1,\dots,y_{\ell-1} \in
B(z,r)$ such that  for all $i \in
\{1,\ldots,\ell-1\}$ we have $\p_n(y_i) \to \xi$ and
\be
|y_i-\xi_i|<|z-\xi_i|.
\ee
We can also find $y_\ell$ such that $\psi_n(y_\ell)\to\xi$ and
$$|y_\ell-\xi_\ell|>|z-\xi_\ell|.$$ Since
for all $i$ the sequence $\p_n(y_i)$ converges to $\xi$ there
exists $N$ such that $\p_n(y_i)=\xi$ for all $i\in \{1,\dots,\ell\}$
and all $n>N$.

There exists $r'>0$ such that for all $y\in B(z,r')$:
\begin{mylist}
\item $|y_i-\xi_i|<|y-\xi_i|$ for all $i<\ell$,

\item $|y_\ell-\xi|>|y-\xi|$, and

\item $|y-\eta|>|y-\xi|$  for all $\eta\in\Xi\setminus
\{\xi_1,\ldots, \xi_\ell\}$.
\end{mylist}

\paragraph{Claim.} For all
$n>N$ and all $y \in B(z,r')$ we have that $\p_n(y)=\xi$ or
$\p_n(y)=\Delta$.

Suppose that the claim does not hold for some $y$ and $n>N$.  If
$\psi_n(y)=\infty$ or if
$\psi_n(y)=\eta\in\Xi\setminus\{\xi_1,\ldots,\xi_\ell\}$ then
$(y,\xi)$ would be unstable by (ii) and (iii) above.  On the other
hand if $\psi_n(y)=\xi_i$ where $i<\ell$ then by (i) and (ii)
$(\xi_i,y_i)$ would be an unstable pair.  Thus the claim is
established.

As for every $n$ the set $\p_n^{-1}(\xi)$ is open and
 $\leb\p_n^{-1}(\Delta)=0$, we deduce from the claim and the
 fact that $\psi$ is a canonical allocation that $\p_n(y)=\xi$ for
 all $y \in B(z,r')$ and for all $n>N$. Thus $\p_n(z)\to \xi.$
\end{pfof}

\section*{Open Problems}

\begin{mylist}

\item {\bf Critical behavior in dimension two and higher.}
What is the tail behavior of
$X$ or $R^*$ for the critical Poisson model? In particular, give
any quantitative upper bound on $\prob(X>r)$ as $r\to\infty$
for $d\geq 2$.

\item {\bf Critical behavior in one dimension.}  Can the critical
model be analyzed exactly in the case $d=1$?  Which moments of $X$ are
finite?  The variant model in which each site is only allowed to be
allocated to a center to its right can be analyzed exactly via of the
function $F$ from Section \ref{sec:1d}.  The method may be found in
\cite{liggett-tagged}, in a slightly different context.  For this
model, $\expe X^\nu<\infty$ if and only if $\nu<1/2$.

\item {\bf Explicit non-critical bounds.}  Give explicit bounds on the
exponential decay rates for the subcritical and supercritical models for
general appetite and dimension.

\item {\bf Supercritical radius.} Does $(R^*)^d$ have exponentially
decaying tail for the supercritical model in dimension $d\geq 2$?

\end{mylist}

\section*{Acknowledgements}
We thank Alan Hammond for suggesting the subcritical and
supercritical models, and for valuable conversations.
Christopher Hoffman and Yuval Peres thank
IMPA in Rio de Janeiro, where some of this work was done.

\bibliography{phd}

\bigskip \noindent
{\bf Christopher Hoffman}: \texttt{hoffman(at)math.washington.edu} \\
Department of Mathematics\\
University of Washington\\
Seattle, WA 98195, USA

\bigskip \noindent
{\bf Alexander E. Holroyd}: \texttt{holroyd(at)math.ubc.ca}\\
Department of Mathematics \\
University of British Columbia \\
Vancouver, BC V6T 1Z2, Canada

\bigskip \noindent
{\bf Yuval Peres}: \texttt{peres(at)stat.berkeley.edu} \\
Departments of Statistics and Mathematics\\
UC Berkeley\\
Berkeley, CA 94720, USA

\end{document}